\documentclass[12pt,a4paper]{amsart}
\usepackage{graphics}
\newtheorem{thm}{Theorem}[section]
\newtheorem{lem}[thm]{Lemma}
\newtheorem{cor}[thm]{Corollary}
\makeatletter
\def\putfig{\@ifnextchar[{\opputfig}{\noopputfig}}
\makeatother
\def\opputfig[#1]#2{\begin{figure}[hbt]
  \centerline{\scalebox{#1}{\includegraphics{#2.eps}}}
  \caption{\space}\label{f:#2}
  \end{figure}}
\def\noopputfig#1{\begin{figure}[hbt]
  \centerline{\includegraphics{#1.eps}}
  \caption{\space}\label{f:#1}
  \end{figure}}
\newcounter{fig}
\newcommand\R{\mathbf{R}}
\allowdisplaybreaks

\begin{document}

\title[Positive presentations of the braid group]
{Positive presentations of the braid groups\\ and the embedding
problem}

\author{Jae Woo Han}
\address{Department of Mathematics\\
Korea Advanced Institute of Science and Technology\\ Taejon,
305--701\\ Korea}
\email{jwhan\char`\@knot.kaist.ac.kr}

\author{Ki Hyoung Ko}
\address{Department of Mathematics\\
Korea Advanced Institute of Science and Technology\\ Taejon,
305--701\\ Korea}
\email{knot\char`\@knot.kaist.ac.kr}

\maketitle

\begin{abstract}
A large class of positive finite presentations of the braid
groups is found and studied. It is shown that no presentations
but known exceptions in this class have the property that
equivalent braid words are also equivalent under positive
relations.
\end{abstract}

\section{Introduction}
We will discover a large class of positive finite presentations of
the $n$-braid group. Roughly speaking, a presentation in this class
has a set of generators determined by a connected graph with $n$ vertices
immersed in a plane in such a way that each pair of edges intersects
at most once. Our class of
presentations includes all known positive finite presentations
such as the Artin presentation\cite{Ar}, the band-generator
presentation in \cite{BKL,KKL} and all of Sergiescu's
presentations in \cite{Se}. Our first objective is to find a
(minimal) collection of positive relations among generators given
by a graph as above so that it becomes a presentation of the $n$-braid groups.

The semigroup of all positive words
plays a crucial role in the
solutions to the word problem and the conjugacy problem for the
$n$-braid group given by
Garside\cite{Ga}, Thurston\cite{Ep}, and Elrifi-Morton\cite{EM}.
Their solutions are based on the property that the semigroup
embeds in the whole group. Birman-Ko-Lee's recent work\cite{BKL}
also requires this embedding property in order to give a fast
solution to the word problem in the band-generator presentation.
Our second objective is to prove that with a few exception all presentations
in our class have the embedding property. The
Artin and the band-generator presentations are the only
nontrivial exceptions. Consequently these two among all
presentations are natural and interesting for further study.

\section{Immersed graphs and known presentations}

The $n$-braid group $B_n$ is the group of isotopy classes of
orientation-preserving automorphisms of an $n$ punctured plane
$\R^2$ that fix the set of punctures and the outside of a disk
containing all punctures. It is customary to place punctures
on the $x$-axis and equally spaced. Ignoring the punctures,
such an automorphism is a homeomorphism of $\R^2$ permuting the
punctures and so is isotopic to the identity map of $\R^2$. A
geometric braid is determined by taking the trace of the punctures
under this isotopy. When the punctures are not in the customary
location, there is a choice of homeomorphisms of $\mathbf R^2$
that send the punctures to their customary location and a
geometric braid is uniquely determined up to conjugation in this
case.

Since the symmetric group over $n$ letters is naturally a quotient
of $B_n$, an automorphism that exchanges two punctures is a
candidate for a generator of $B_n$ and a typical automorphism of
this type can be easily depicted as an arc connecting the two
punctures so the automorphism exchanges two punctures clockwise or
counterclockwise inside an annulus lying along the arc. See Figure
\ref{f:string} for an example. This type of braids will be called
the {\em half twist along the given arc}.

\putfig{string}

Therefore by regarding the punctures as vertices and a collection
of these arcs as edges, we have a graph $\Gamma$ immersed in
$\R^2$ satisfying the following properties:

\begin{enumerate}
\item There are no loops;
\item Edges have no self-intersections;
\item Two distinct edges intersect either at interior points transversely
or at common end points.
\item There is at least a vertex inside a ``pseudo digon",  that is,
a region cobounded by two subarcs from two edges so that corners of the
region are either two interior intersection points or one vertex
and one interior intersection point.
\end{enumerate}

If the number of intersections among edges is minimized,
for example, by making all edges geodesics in the hyperbolic
structure of the punctured plane, the condition (4) holds automatically.
Two immersed graphs are regarded as equivalent if one graph can be
transformed to the other by an orientation-preserving
homeomorphism of $\mathbf R^2$. Thus we may assume that
vertices lie on the customary locations if necessary.

Throughout this article, graphs mean immersed graphs in $\mathbf
R^2$ satisfying the condition discussed. If graphs are embedded in
$\mathbf R^2$, we call them {\em planar}. The set of edges in a
graph $\Gamma$ with $n$ vertices corresponds to a set of elements
in $B_n$ as in Figure~\ref{f:string}. Two equivalent graphs
determine two sets of braids that differ by an inner automorphism
of $B_n$. If the set corresponding to edges generates $B_n$, we
say that $\Gamma$ {\em generates} $B_n$ and we are interested in
graphs that generate $B_n$. Since an edge of $\Gamma$ corresponds
to a transposition in the symmetric group, a graph that generates
$B_n$ must be connected in order to permute any two vertices. But
there are connected graphs that do not generate $B_n$. The problem
of deciding when an immersed graph generates the braid group is
considered elsewhere\cite{HK}.

We now discuss how a braid word corresponds to an edge of a graph
with $n$ vertices in the customary location.  The Artin
generator $\sigma_i$(or $\sigma_i^{-1}$) for $i=1,\ldots,n-1$
corresponds to the clockwise(counterclockwise, respectively) half
twist along the straight edge joining the $i$-th and the
$(i+1)$-st vertices. Given an edge of a graph, we simplify the
edge via a sequence of half twists $\sigma_i$ or $\sigma_i^{-1}$
for $i=1,\ldots,n-1$ until it becomes a straight edge
corresponding to, say, $\sigma_k$. The required sequence can be
expressed as a word $W$ in Artin generators in which a half
twist applied later is written on the left. Then the edge that we
started with corresponds to the braid word $W^{-1}\sigma_k W$. The
word $W$ is not uniquely expressed but $W$ as an element of $B_n$
is well-defined. The inverse of this correspondence is similar.
Given a conjugate $W^{-1}\sigma_k W$, apply a sequence of half
twists determined from right to left by $W^{-1}$ to the straight
edge $\sigma_k$. Consequently each edge of a graph with $n$
vertices is uniquely represented by an element in $B_n$ that can
be written as a conjugate of an Artin generator. The following
lemma summarizes this discussion.

\begin{lem}
A braid $\beta$ in $B_n$ is the half twist along an arc joining
two punctures in the plane if and only if $\beta$ can be written
as a conjugate of an Artin generator or its inverse.
\end{lem}

From now on as long as no confusion arises, we will not
distinguish three concepts, namely an edge of a graph, the half
twist along the edge, and a conjugate word expressing the edge.
Since any two Artin generators are conjugate each other, any
half twist along an arc corresponds to a conjugate of a fixed
Artin generator. The following is immediate from the above lemma.

\begin{cor}\label{conjugate}
A graph $\Gamma$ of $n$ vertices generates $B_n$ if and only if
each Artin generator can be expressed as $W\alpha W^{-1}$ for
some edge $\alpha$ and some word $W$ on edges of $\Gamma$.
\end{cor}

In the view of the above corollary, it is important to know how
the conjugate of an edge by another edge in a graph looks like.
Let $\alpha$ and $\beta$ be words in $B_n$ that express edges of a
graph. If the edges $\alpha$ and $\beta$ do not intersect each
other and are not adjacent, they commute and so
$\alpha\beta\alpha^{-1}=\beta=\alpha^{-1}\beta\alpha$. If the
edges $\alpha$ and $\beta$ do not intersect each other and are
adjacent, the three edges $\alpha$, $\beta$, and $\alpha\beta\alpha^{-1}
=\beta^{-1}\alpha\beta$ form a triangle counterclockwise and the
three edges $\alpha$, $\beta$, and $\alpha^{-1}\beta\alpha
=\beta\alpha\beta^{-1}$ form a triangle clockwise as in
Figure~\ref{f:triangle}. We note that there should be no other
vertex inside the triangles and this can be achieved by drawing
thin triangles inside a sufficiently small neighborhood of the
union of two edges $\alpha$ and $\beta$. When the edges $\alpha$
and $\beta$ intersect at an interior point, one can still describe
the edges $\alpha\beta\alpha^{-1}$ and $\alpha^{-1}\beta\alpha$
but we do not need them in this article.

\putfig{triangle}

Give a set $E$ of generators of the $n$-braid group $B_n$, a {\em
positive word} in $E$ is a product of positive powers of
generators in $E$. A presentation of a group is {\em finite} if
there are finitely many generators and relations. A presentation
is {\em positive} if all defining relations are equations of
positive words in generators.

We now introduce the graphs that corresponds to some of known
positive finite presentations of the braid groups.
The set of Artin generators
$\sigma_1,\sigma_2,\sigma_3,\ldots,\sigma_{n-1}$ of $B_n$ forms
the graph in Figure~\ref{f:artin} and a minimal set of defining
relations is given by
\begin{eqnarray*}
&\sigma_i\sigma_j\sigma_i = \sigma_j\sigma_i\sigma_j \qquad
&\mbox{if}\quad  |i-j| = 1\\ &\sigma_i\sigma_j = \sigma_j\sigma_i
\qquad &\mbox{if}\quad |i-j| > 1.
\end{eqnarray*}

\putfig{artin}

\noindent The graphs in Figure~\ref{f:artin} will be called the {\em
Artin} graph.

Sergiescu \cite{Se} showed how a finite presentation of
$B_n$ can be obtained from any planar connected graph with $n$
vertices. He described a sufficient set of positive relations that
depend only on the geometry of a given planar graph. We will give a
minimal set of positive relations as a corollary of the main
theorem in \S 3.
Let $n$ be the number of vertices in $\Gamma$. By choosing a fixed
$n$ points in $\mathbf{R}^2$ and a homeomorphism of $\mathbf{R}^2$
that sends vertices of $\Gamma$ to the $n$ fixed points, the group
$B_\Gamma$ introduced in \cite{Se} is identified with the $n$-braid
group $B_n$. These presentations of $B_n$ will be called {\em Sergiescu's
presentations.} The relations in a Sergiescu's presentation are
highly redundant but they are useful when we need to find a
relation locally given by a planar graph.

Recently a new presentation of $B_n$ called the band-generator
presentation has been developed by Birman-Ko-Lee\cite{BKL}. This
presentation has $n\choose 2$ generators $a_{ts}$ for $1\le s<t\le
n$ corresponds to the graph in Figure~\ref{f:band} and defining
relations:
\begin{eqnarray*}
&a_{ts}a_{sr}=a_{tr}a_{ts}=a_{sr}a_{tr} \quad &{\rm for \  all} \
t,s,r
 {\rm \  with} \ \ n\geq t>s>r \geq 1\\
&a_{ts}a_{rq}=a_{rq}a_{ts}\quad &\mbox{\rm if } \
(t-r)(t-q)(s-r)(s-q)>0.
\end{eqnarray*}

The Artin generator and the band-generators are related as
$$a_{ts}=(\sigma_{t-1}\sigma_{t-2}\cdots\sigma_{s+1})\sigma_s
(\sigma_{s+1}^{-1}\cdots\sigma_{t-2}^{-1}\sigma_{t-1}^{-1})$$ and
$$\sigma_t=a_{(t+1)t}.$$ The set of band-generators are depicted as
the graph in Figure~\ref{f:band} that will be called the
{\em inner complete} graph.

\putfig{band}

The relations $a_{ts}a_{sr}=a_{tr}a_{ts}=a_{sr}a_{tr}$ will
be called a {\em triangular} relation. A triangular relation is
derived whenever a new generator is introduced by means of a
conjugation of an edge by an adjacent edge. Triangular relations
serve as building blocks of positive relations in the braid groups.

\section{Positive Presentations from linearly spanned graphs}

After an immersed graph is turned into a planar graph by regarding
all interior intersections as vertices, a region bounded by a closed
edge-path is called a {\em pseudo face} if the edge-path contains
at least one vertex that is an interior intersection of two edges
as in Figure~\ref{f:pseudo}.

\putfig{pseudo}

A graph $\Gamma$ is said to be {\em linearly spanned} if it is
connected and  there is no vertex in any pseudo face of
$\Gamma$. A connected subgraph of a linearly spanned graph is
clearly linearly spanned. Two edges of a linearly spanned graph
intersect each other at most once since any pseudo digon cobounded
by two edges can not exist in a linearly spanned graph. All planar
graphs and all subgraphs of the inner-complete graph are linearly
spanned. For example the graph on the left in
Figure~\ref{f:linear} is neither planar nor a subgraph of the
inner-complete graph but it is linearly spanned, and the graph on
the right is not linearly spanned

\putfig{linear}

The following lemma justifies the terminology "linearly spanned".
\begin{lem}
A linearly spanned graph that is a tree is equivalent to a
subgraph of the inner-complete graph.
\end{lem}
\begin{proof}
Choose a point $x$ far away from the graph. Since the graph is a tree
and no vertices are surrounded by edges, there are  $n$ arcs that
join $x$ to each vertex and are disjoint each other and are
disjoint from the graph as in Figure~\ref{f:tree}. Choose a new
horizontal axis disjoint from the graph and move each vertex along
each arc by a homeomorphism so that it lies on the new axis as
Figure~\ref{f:tree}. Then the result is a subgraph of the inner
complete graph.

\putfig{tree}
\end{proof}

\begin{cor}
A linearly spanned graph $\Gamma$ with $n$ vertices generates
$B_n$.
\end{cor}
\begin{proof}
A maximal tree of $\Gamma$ is equivalent to a connected subgraph
of the inner-complete graph from which we can obtain the inner
complete graph by adding missing edges that are conjugates of an
existing edge by an adjacent edge.
\end{proof}

\begin{lem}\label{tree}
Let $T$ be a linearly spanned tree with $n$ vertices. Then $T$
generates $B_n$ with $\frac{(n-1)(n-2)}2$ positive relations
described as in the proof below.
\end{lem}
\begin{proof}
Induction on $n$. It is trivial when $n=2$. Suppose a linear
spanned tree $T'$ with $n-1$ vertices generates $B_{n-1}$ with
$\frac{(n-2)(n-3)}2$ positive relations. We add a new vertex and a
new edge $\alpha$ to $T'$. For each edge $\beta$ of $T'$, we will
have a new positive relation so that $n-2$ new positive relations
will be added. In the following we use Tietze transformations
\cite{Ti} that add and delete a generator(s) denoted by $\lambda$
or $\mu$ to utilize triangular relations.
\begin{enumerate}
\item If $\alpha$ and $\beta$ have no intersection, add the
positive relation $$\alpha\beta=\beta\alpha.$$

\item For each set of edges $\beta_1,\ldots,\beta_m$ simultaneously
adjacent to $\alpha$ at a vertex of valency $m+1$ as in
Figure~\ref{f:induct1}(a), add $m$ positive relations
%$$\alpha\beta_m\alpha = \beta_m\alpha\beta_m,
%\alpha\beta_1\beta_m\alpha = \beta_m\alpha\beta_1\beta_m, \ldots,
%\alpha\beta_{m-1}\beta_m\alpha = \beta_m\alpha\beta_{m-1}\beta_m
%$$
\begin{eqnarray*}
\alpha\beta_m\alpha &=& \beta_m\alpha\beta_m\\
\alpha\beta_1\beta_m\alpha &=& \beta_m\alpha\beta_1\beta_m\\
&\vdots& \\
\alpha\beta_{m-1}\beta_m\alpha &=&
\beta_m\alpha\beta_{m-1}\beta_m
\end{eqnarray*}
\item If $\alpha$ and $\beta$ intersect and form a pseudo face as in
Figure~\ref{f:induct1}(b), add the positive relation
\begin{eqnarray*}
\lefteqn{\beta_1 \beta_2 \cdots \beta_m \alpha \beta \beta_1^2
\beta_2 \cdots \beta_m \alpha \beta \beta_1 \beta_2 \cdots
\beta_{m-1}}\\ &=& \beta_2 \beta_3 \cdots \beta_m \alpha \beta
\beta_1^2 \beta_2 \cdots \beta_m \alpha \beta \beta_1 \beta_2
\cdots \beta_m
\end{eqnarray*}
that is derived from $\beta \lambda \beta = \lambda \beta \lambda$
where
\begin{eqnarray*}
\lambda &=& \beta_1 \beta_2 \cdots \beta_m \alpha
            \beta_m^{-1} \cdots \beta_2^{-1} \beta_1^{-1}\\
        &=& \beta^{-1} \beta_m^{-1} \cdots \beta_2^{-1}
            \beta_1 \beta_2 \cdots \beta_m \alpha
\end{eqnarray*}

\putfig{induct1}

\item If $\alpha$ and $\beta$ intersect and form a pseudo face as in
Figure~\ref{f:induct2}(a), add the positive relation $$ \beta_1
\beta_2 \cdots \beta_m \alpha \beta \beta_1 \beta_2 \cdots
\beta_{m} = \beta_2 \beta_3 \cdots \beta_m \alpha \beta \beta_1
\beta_2 \cdots \beta_m \alpha$$ that is derived from $\beta
\lambda= \lambda \beta$ where
\begin{eqnarray*}
\lambda &=& \beta_1 \beta_2 \cdots \beta_m \alpha
            \beta_m^{-1} \cdots \beta_2^{-1} \beta_1^{-1}\\
        &=& \beta^{-1} \beta_m^{-1} \cdots \beta_2^{-1}
            \beta_1 \beta_2 \cdots \beta_m \alpha
\end{eqnarray*}

\item If $\alpha$ and $\beta$ intersect and form a pseudo face as in
Figure~\ref{f:induct2}(b),  add the positive relation
\begin{eqnarray*}
\lefteqn{\gamma_2 \cdots \gamma_\ell \alpha \beta \beta_1 \cdots
\beta_m
 \gamma_1 \cdots \gamma_\ell \alpha \beta \beta_1 \cdots \beta_{m-1}}\\ &=&
\beta_1 \cdots \beta_m \gamma_1 \cdots \gamma_\ell \alpha \beta
 \beta_1 \cdots \beta_m \gamma_1 \cdots \gamma_\ell
\end{eqnarray*}
that is derived from $\lambda \mu = \mu \lambda$ where
\begin{eqnarray*}
\lambda &=& \beta \beta_1 \beta_2 \cdots \beta_m
            \beta_{m-1}^{-1} \cdots \beta_2^{-1} \beta_1^{-1} \beta^{-1}\\
        &=& \beta_m^{-1} \cdots \beta_1^{-1}
            \beta \beta_1 \beta_2 \cdots \beta_m\\
\mu &=& \gamma_1 \gamma_2 \cdots \gamma_\ell \alpha
            \gamma_\ell^{-1} \cdots \gamma_2^{-1} \gamma_1^{-1}\\
        &=& \alpha^{-1} \gamma_\ell^{-1} \cdots \gamma_2^{-1}
            \gamma_1 \gamma_2 \cdots \gamma_\ell \alpha
\end{eqnarray*}

\end{enumerate}

\putfig{induct2}

\end{proof}

The following is the main theorem of this section and the known
presentation mentioned in the previous section can be obtained
from this.

\begin{thm}\label{circuit}
Let $\Gamma$ be a linearly spanned graph with $n$ vertices. Then
$\Gamma$ generates $B_n$ with $\frac{(n-1)(n-2)}2+k$ positive
relations described as in the proof below where $\Gamma$ has
$n+k-1$ edges.
\end{thm}
\begin{proof}
Choose a spanning tree $T$ of $\Gamma$. Then a linearly spanned
tree $T$ generates $B_n$ with $\frac{(n-1)(n-2)}2$ positive
relations as in Lemma~\ref{tree}. When each edge $\alpha$ in
$\Gamma - T$ are added, a circuit is formed and the circuit give a
new positive relation as follows:

\begin{enumerate}
\item If $\alpha$ forms a circuit with no intersection with other
edges as in Figure~\ref{f:cutopen}, add the positive relation
$$\alpha\beta_1\beta_2 \cdots \beta_{l-1}= \beta_1\beta_2 \cdots
\beta_\ell$$ where we regard that the circuit bounds a polygonal disk
by cutting open all edges inside the circuit as in
Figure~\ref{f:cutopen}.

\putfig{cutopen}

\item If $\alpha$ forms a circuit with some intersections with other
edges as in Figure~\ref{f:face}, add the positive relation
\begin{eqnarray*}
\lefteqn{\beta_{11}\beta_{12} \cdots \beta_{1l_1} \beta_{2l_2}
\alpha
   \beta_{2(l_2-1)} \cdots \beta_{22}} \\ &=&
\beta_{12} \cdots \beta_{1l_1} \beta_{2l_2} \alpha
   \beta_{2(l_2-1)} \cdots \beta_{22} \beta_{21}
\end{eqnarray*}

\putfig{face}
\end{enumerate}

\end{proof}

\section{Embedding problem}
Two positive words $U$, $V$ in a positive presentation will be
said to be {\em positively equivalent} if they are identically
equal or they can be transformed into each other through a
sequence of positive words such that each word of the sequence is
obtained from the preceding one by a single direct application of
the defining relations. And we will write $U \doteq V$ if $U$ and
$V$ are positively equivalent. Given a positive finite
presentation $\langle X\,|\, R\rangle$ of $B_n$, let $B^+_n$ be
the free semigroup generated by $X$ modulo $R$. If any two
equivalent positive words $U$ and $V$ are positively equivalent,
then we say that the semigroup $B^+_n$ {\em embeds} in $B_n$ or
the presentation $\langle X\,|\, R\rangle$ has the {\em embedding
property}. A set $X$ of generators is said to have the {\em
embedding property} if a presentation $\langle X\,|\, R\rangle$
has the embedding property for some finite set $R$ of positive
relations. The {\em embedding problem} of a graph that generates
$B_n$ with positive relations is to decide whether the set of
generators given by the graph has the embedding property. Thus a
graph does not have the embedding property if and only if no
finite set of positive relations over the set of generators given
by the graph form a presentation with the embedding property.

Garside\cite{Ga} showed that the Artin presentation has the
embedding property and Birman-Ko-Lee\cite{BKL,KKL} showed that the
band-generator presentation has the embedding property. We will
show that linearly spanned graphs with more than 3 vertices do not
have the embedding property except these two presentations. There
are 4 possible linearly spanned graphs with 3 vertices as in
Figure~\ref{f:3vertex}.

\putfig{3vertex}

The first two graphs give the Artin and the band-generator
presentations of $B_3$. The last two graphs have multiple edges.
We do not know whether these graphs have the embedding property.
In particular the following presentation with 4 defining
relations: $$\langle\alpha,\beta,\gamma\,|\,\alpha \beta \alpha =
\beta \alpha \beta, \gamma \beta \gamma = \beta \gamma \beta,
\beta^2 \gamma = \alpha \beta^2,\beta \gamma \alpha = \gamma
\alpha \beta\rangle$$ may have the embedding property. But all
known techniques as in \cite{Bi,BKL,CP,KKL,Re} fail to apply to
this example. We will try to avoid these unknown exceptions in the
following discussion by allowing no multiple edges.

A subgraph $\Gamma'$ of a graph $\Gamma$ is said to be {\em full}
if every edge in $\Gamma$ joining two
vertices in $\Gamma'$ is also in $\Gamma'$.

\begin{thm}\label{circle}
Let $\Gamma'$ be a connected full subgraph of a linearly spanned graph
$\Gamma$. If a graph $\Gamma$ has the embedding property and there
is a circle $C$ such that $C$ contains $\Gamma'$ inside and all
vertices in $\Gamma-\Gamma'$ lie outside, then $\Gamma'$ also has
the embedding property.
\end{thm}
\begin{proof}
Let $X$ and $X'$ be the set of generators given by $\Gamma$ and
$\Gamma'$, respectively and let $R$ be a finite set of positive
relations on $X$ such that $\langle X\,|\, R\rangle$ has the
embedding property.
Let $R'$ be the set of relations on $X'$ which is
a ``full" subset of $R$ in the sense that $R'$ contains all
relations in $R$ written on $X'$. We will prove by contradiction that
$\langle X'\,|\, R'\rangle$ has the embedding property.
Suppose that there exists a pair of
positive words $U,V$ on $X'$ such that $U, V$ are equivalent
in the braid group but are not positively equivalent under $R'$.
Since $\langle X\vert R\rangle$ has the
embedding property, $U \doteq V$ under $R$. Thus there is a
sequence of positive words $W_1,\ldots,W_k$ over $X$ such that
$U\doteq W_1\doteq\cdots\doteq W_k\doteq V$ under $R$ and each
positive equivalence is obtained by one direct application of
relations in $R$. The sequence must contain at least one positive word,
say $W_i$, that is not written solely on $X'$, otherwise $U \doteq V$ under $R'$.
In the view of Lemma~\ref{tree} and Theorem~\ref{circuit}, it is
impossible that $R$ contains any relation $W=W'$ such that $W'$ is a positive word
over $X'$ and $W$ is a positive word over the edges that are not
incident to any vertex of $\Gamma'$. Furthermore $\Gamma'$ is a
full subgraph of $\Gamma$. Thus $W_i$ must contain an edge that
joins a vertex $v$ in $\Gamma'$ to a vertex $w$ not in $\Gamma'$.
Let $\beta$ be such an occurrence that comes last in $W_i$.
Since $\Gamma'$ is connected, the vertex $v$ is joined to another vertex
$u$ in $\Gamma'$ by an edge $\alpha$. As an automorphism of the
punctured plane, $U$ does not change the circle $C$ because the edges
for $U$ never touch $C$. On the other hand we will show $W_i$ must
change $C$ and this is a contradiction because $U$ and $W_i$ are
isotopic as automorphisms of the punctured disk.
The automorphism $W_i$ is the composition of the counterclockwise half
twists along edges in $W_i$. Then the counterclockwise half twist
along $\beta$ creates an intersection $x$ of $C$ with the edge
$\alpha$. The intersection $x$ can disappear only via a clockwise
half twist along an edge incident at either $u$ or $v$ as in
Figure~\ref{f:circle}. But all of half twists in $W_i$ are
counterclockwise because $W_i$ is a positive word over edges,

\putfig{circle}
\end{proof}

We think the above theorem also holds when we replace the
condition "linearly spanned" by "connected".

\begin{lem}\label{oneinter}
A linearly spanned graph $\Gamma$ with 4 vertices has at most one
intersection among its edges.
\end{lem}
\begin{proof} An intersection between two adjacent edges always
creates a pseudo face that must contain a vertex and so there is
no intersection between adjacent edges in a linearly spanned
graph. Suppose $\Gamma$ has more than one intersection.
Figure~\ref{f:oneinter} shows a typical situation with one
intersection and another intersection $x$. The edge $\alpha$ must
join vertices $v_1$ and $v_3$, otherwise $x$ is an intersection
between two adjacent edges. One can easily check that all
possibilities of completing $\alpha$ create a pseudo face
containing at least a vertex.

\putfig{oneinter}
\end{proof}

%\begin{lem}\label{length}
%Let $\Gamma$ be a linearly spanned graph. Suppose that for each
%$k=1,2,\ldots$, there are positive words $W,W'$ over edges in
%$\Gamma$ such that $W=W'$, the word length $|W|=|W'|=k+c$ for some
%constant $c$, and no proper subwords of $W$ is equivalent to a
%distinct positive word over edges in $\Gamma$. Then the graph
%$\Gamma$ does not have the embedding property.
%\end{lem}
%\begin{proof}
%Let $X$ be the set of generators given by $\Gamma$ and $\langle
%X\,|\,R\rangle$ be any positive finite presentation of the braid
%group. Choose a large $k$ such that $W=W'$ is longer than any
%relation in $R$. Since any proper subword of $W$ is unique as a
%positive word over $X$, $W$ can not be positively equivalent to
%$W'$ unless $W=W'$ itself is contained in $R$, which is impossible
%\end{proof}

\begin{lem}\label{length}
Let $\Gamma$ be a linearly spanned graph and $X$ be the set of
generators given by edges of $\Gamma$. Suppose that for each
$k=1,2,\ldots$, there are positive words $W\equiv\alpha V\beta$ and
$W'\equiv\alpha' V'\beta'$ over $X$ for
$\alpha,\alpha',\beta,\beta'$ in $X$ such that
\begin{enumerate}
\item[(i)] $W=W'$;
\item[(ii)] the word length $|W|=|W'|=k+c$ for some constant $c$;
\item[(iii)] $\alpha V\not=\alpha'P$ and $V\beta \not=Q\beta'$
for any positive words $P,Q$ over $X$
\end{enumerate}
Then the graph $\Gamma$ does not have the embedding property.
\end{lem}
\begin{proof}
Let $\langle X\,|\,R\rangle$ be any positive finite presentation
of the braid group. The hypothesis (iii) implies that any shorter
positive relation than $W=W'$ itself can not make $W$ positively
equivalent to $W'$. Choose a large $k$ such that $W=W'$ is longer
than any relation in $R$. Then $W$ is not positively equivalent to
$W'$ over $R$.
\end{proof}

The following theorem completely determines when a linearly spanned
graph with 4 vertices and no multiple edges has the embedding
property. In the proof given below, $P(\alpha_1,\ldots,\alpha_k)$
or $Q(\alpha_1,\ldots,\alpha_k)$ will denote a positive word over
the generators $\alpha_1,\ldots,\alpha_k$.

\begin{thm}\label{4vertex}
Among all linearly spanned graphs with 4 vertices and no multiple
edges, only two of them, the Artin graph and the inner complete
graph as in Figure~\ref{f:embed}, have the embedding property.

\putfig{embed}
\end{thm}
\begin{proof}
The graphs in Figure~\ref{f:embed} give the Artin
presentation and the band-generator presentation. So they have the
embedding property.

In order to show that all other linearly spanned graphs with 4
vertices do not have the embedding property, we appeal to
Lemma~\ref{length}. To check the hypothesis (iii) of the lemma, we
use the fact that the band-generator presentation has the
embedding property so that the positive equivalence in the
presentation is the same as the equivalence in the braid group. We
also utilize the left and right cancellation theorem and the left
and right canonical forms in the band-generator presentation in
\cite{BKL,KKL}.

In the view of Lemma~\ref{oneinter}, the graphs being considered
are divided into two types: planar graphs and graphs with one
intersection among edges.

\subsection*{Planar graphs}

Planar graphs with 4 vertices can be divided further into three types:
graphs containing  neither a triangle nor a rectangle, graphs
containing at least a rectangle, and graphs containing at least a
triangle but no rectangle, where a triangle (or a rectangle) is
non-degenerate, that is, must have 3 (or 4, respectively) vertices
and must contain no other vertices inside.

\subsubsection*{\bf I. Graphs containing  neither a triangle nor a
rectangle}

This type is further divided into two types (i) and (ii).

\paragraph{\bf (i) Graphs containing a vertex adjacent to all of remaining
three vertices}

\putfig{star1}

We may have two possible graphs as in Figure~\ref{f:star1} up to
equivalence. For each of these two graphs, let $W=\alpha_1
\alpha_2 \alpha_3^k \alpha_1$ and $W'=\alpha_2 \alpha_3^k \alpha_1
\alpha_2$. We will show that $W$ and $W'$ satisfy
Lemma~\ref{length}. Add the edges $\lambda, \mu, \nu$ so that
these edges together with $\alpha_1, \alpha_2, \alpha_3$ form an
inner-complete graph. Then $$\alpha_1 \alpha_2 \alpha_3^k \alpha_1
=\alpha_2 \lambda \alpha_3^k \alpha_1
=\alpha_2\alpha_3^k\lambda\alpha_1
=\alpha_2\alpha_3^k\alpha_1\alpha_2$$ and so $W=W'$.

Suppose $\alpha_1 \alpha_2 \alpha_3^k=\alpha_2P$ for some positive
word $P$ in the band-generator presentation. Then $U\doteq
\lambda\alpha_3^k$ by the left cancellation in the band-generator
presentation. But $\lambda\alpha_3^k$ may only start with
$\lambda$ or $\alpha_3$ which commute. Thus $P$ can not be
written over $\alpha_1, \alpha_2, \alpha_3$. Similarly $\alpha_2
\alpha_3^k \alpha_1\not=Q(\alpha_1, \alpha_2, \alpha_3)\alpha_2$.
Thus $W$ satisfies (iii) of Lemma~\ref{length}.

Due to the triangular relation $\lambda\mu=\mu\alpha_4
=\alpha_4\lambda$, $W$ can not be equivalent to a positive word
over $\alpha_1, \alpha_2, \alpha_3, \alpha_4$ that contains
$\alpha_4$ since $W$ is not positively equivalent to a positive
word in the band-generator presentation that contains the subword
$\lambda\mu$. Thus $W,W'$ also satisfy the hypothesis of
Lemma~\ref{length} over $\alpha_1, \alpha_2, \alpha_3, \alpha_4$.

\paragraph{\bf (ii) Graphs containing no vertex adjacent to all of remaining
three vertices} All possible graphs except the Artin graph have
multiple edges as in Figure~\ref{f:except}.

\putfig{except}

\subsubsection*{\bf II. Graphs containing at least a rectangle}
We may have two possible graphs as in Figure~\ref{f:cycle} up to
equivalence. For each of these five graphs, let $W=\alpha_1
\alpha_2^k \alpha_3$ and $W'=\alpha_3 \alpha_4^k \alpha_1$.

\putfig{cycle}

Add the edges $\lambda, \mu$ so that these edges together with
$\alpha_1, \alpha_2, \alpha_3, \alpha_4$ form an inner-complete
graph. Then $$\alpha_1 \alpha_2^k \alpha_3 = \alpha_1
\alpha_2^{k-1} \alpha_3\lambda = \alpha_1\alpha_3\lambda^k
=\alpha_3\alpha_1\lambda^k =\alpha_3\alpha_4\alpha_1\lambda^{k-1}
=\alpha_3 \alpha_4^k \alpha_1$$ and so $W=W'$.

Suppose $\alpha_1 \alpha_2^k=\alpha_2P$ for some positive word $P$
in the band-generator presentation. Then $P\doteq
\mu\alpha_2^{k-1}$ by the left cancellation in the band-generator
presentation. Thus $P$ can not be written over $\alpha_1,
\alpha_2, \alpha_3, \alpha_4$. Similarly $\alpha_2^k \alpha_3
\not=Q(\alpha_1,\ldots, \alpha_4)\alpha_1$. Thus $W$ satisfies
(iii) of Lemma~\ref{length}.

Due to the triangular relation $\alpha_1\alpha_4=\alpha_4\alpha_5
=\alpha_5\alpha_1$, $W$ can not be equivalent to a positive word
over $\alpha_1,\ldots, \alpha_5$ that contains $\alpha_5$ since
$W$ is not positively equivalent to a positive word in the
band-generator presentation that contains the subword
$\alpha_1\alpha_4$. Thus $W,W'$ also satisfy the hypothesis of
Lemma~\ref{length} over $\alpha_1,\ldots, \alpha_5$.

\subsubsection*{\bf III. Graphs containing at least a triangle but
no rectangle} In this case, graphs are divided two types by the
number of triangles.
\paragraph{\bf (i) Graphs containing one triangle}
We have one possible graph as in Figure~\ref{f:outer1} up to
equivalence. Let
$W=\alpha_1\alpha_2\alpha_3\alpha_4\alpha_1\alpha_2^k\alpha_3\alpha_4$
and
$W'=\alpha_2\alpha_3^k\alpha_4\alpha_1\alpha_2\alpha_3\alpha_4\alpha_1$.

\putfig{outer1}

Add the edges $\lambda, \mu$ so that these edges together with
$\alpha_1, \alpha_2, \alpha_3, \alpha_4$ form an inner-complete
graph. Then
\begin{eqnarray*}
\alpha_1\alpha_2\alpha_3\alpha_4\alpha_1\alpha_2^k\alpha_3\alpha_4
&=& \alpha_1\alpha_2\alpha_3\alpha_1\mu\alpha_2^k\alpha_3\alpha_4
= \alpha_1\alpha_2\alpha_1\alpha_3\mu\alpha_2^k\alpha_3\alpha_4\\
&=& \alpha_2\alpha_1\alpha_2\alpha_3\mu\alpha_2^k\alpha_3\alpha_4
= \alpha_2\alpha_1\alpha_2\alpha_3\alpha_2^k\mu\alpha_3\alpha_4 \\
&=& \alpha_2\alpha_1\alpha_3^k\alpha_2\alpha_3\mu\alpha_3\alpha_4
= \alpha_2\alpha_1\alpha_3^k\alpha_2\mu\alpha_3\mu\alpha_4 \\ &=&
\alpha_2\alpha_1\alpha_3^k\mu\alpha_2\alpha_3\mu\alpha_4 ~ =
\alpha_2\alpha_1\alpha_3^k\mu\alpha_2\alpha_3\alpha_4\alpha_1 \\
&=& \alpha_2\alpha_3^k\alpha_1\mu\alpha_2\alpha_3\alpha_4\alpha_1
=
\alpha_2\alpha_3^k\alpha_4\alpha_1\alpha_2\alpha_3\alpha_4\alpha_1
\end{eqnarray*}
and so $W=W'$.

Let $\alpha_1\alpha_2\alpha_3\alpha_4\alpha_1\alpha_2^k\alpha_3 =
\alpha_2P$. Then $P =
\alpha_1\alpha_2\alpha_3\mu\alpha_2^k\alpha_3\doteq
\alpha_1\alpha_2\alpha_3\alpha_2^k\alpha_3\lambda$. The right
canonical form of $\alpha_1\alpha_2\alpha_3\alpha_2^k\alpha_3$ is
$\alpha_3^{k-3}(\alpha_1\alpha_3)\alpha_4(\alpha_3\alpha_2)(\alpha_3\alpha_2)$.
So $\alpha_1\alpha_2\alpha_3\alpha_2^k\alpha_3$ may finish with
$\alpha_2, \alpha_3, \alpha_4$ over the band-generators and so
$\lambda$ or $\mu$ must appear in $P$. Similarly
$\alpha_2\alpha_3\alpha_4\alpha_1\alpha_2^k\alpha_3\alpha_4\not=
Q(\alpha_1,\ldots,\alpha_4)\alpha_1$.

\paragraph{\bf (ii) Graphs containing two triangles}
In this case triangles must share at least 1 edge since we
consider only 4 vertices. However if two triangles share 2 edges
then it must contain a multiple edge. Thus we have two possible
graphs as in Figure~\ref{f:trian1}. For these two graphs, let
 $W=\alpha_2\alpha_3^k\alpha_5\alpha_2\alpha_3$ and
 $W'=\alpha_4\alpha_1^k\alpha_5\alpha_4\alpha_1$.

\putfig{trian1}

Add the edges $\lambda$ so that these edges together with
$\alpha_1,\ldots, \alpha_5$ form an inner-complete graph. Then
\begin{eqnarray*}
\alpha_2\alpha_3^k\alpha_5\alpha_2\alpha_3 &=&
\lambda^k\alpha_2\alpha_5\alpha_2\alpha_3 =
\lambda^k\alpha_2\alpha_5\alpha_3\lambda =
\lambda^k\alpha_2\alpha_4\alpha_5\lambda \\ &=&
\lambda^k\alpha_4\alpha_2\alpha_5\lambda =
\alpha_4\alpha_1^k\alpha_2\alpha_5\lambda =
\alpha_4\alpha_1^k\alpha_5\alpha_1\lambda \\ &=&
\alpha_4\alpha_1^k\alpha_5\alpha_4\alpha_1
\end{eqnarray*}
and so $W=W'$.

The left canonical form of $\alpha_2 \alpha_3^k \alpha_5 \alpha_2$
is $(\alpha_2 \alpha_3)\alpha_3^{k-1}\alpha_5 \alpha_2$ and so it
can not start with $\alpha_4$ over the band-generators. Thus
$\alpha_2 \alpha_3^k \alpha_5
\alpha_2\not=\alpha_4P(\alpha_1,\ldots, \alpha_5)$. The right
canonical form of $\alpha_3^k\alpha_5\alpha_2\alpha_3$ is
$\alpha_3^k\alpha_5(\alpha_2\alpha_3)$ and so it can not finish
with $\alpha_1$. Thus $\alpha_3^k \alpha_5 \alpha_2\alpha_3
\not=Q(\alpha_1,\ldots, \alpha_5)\alpha_1$.

Due to the triangular relation $\alpha_1\alpha_4=\alpha_4\alpha_6
=\alpha_6\alpha_1$, $W$ can not be equivalent to a positive word
over $\alpha_1,\ldots, \alpha_5$ that contains $\alpha_6$ since
$W$ is not positively equivalent to a positive word in the
band-generator presentation that contains the subword
$\alpha_1\alpha_4$. Thus $W,W'$ also satisfy the hypothesis of
Lemma~\ref{length} over $\alpha_1,\ldots, \alpha_5$.

\subsection*{Graphs with one intersection}

The edges $\alpha_1$ and $\alpha_3$ that intersect each other are
transformed to the diagonals of a rectangle. Then 4 more edges are
needed to make an inner-complete graph. We have four
distinct types of graphs, depending on how many edges are missing
from the inner-complete graph.

\subsubsection*{\bf I. Three edges are missing}

Figure~\ref{f:nonplane} is the only possible graph in this case.
Let $W=\alpha_3\alpha_1\alpha_2^k\alpha_3\alpha_1\alpha_2$ and
$W'=\alpha_2\alpha_3\alpha_1\alpha_2^k\alpha_3\alpha_1$.

 \putfig{nonplane}

Add the edges $\lambda,\mu,\nu$ so that these edges together with
$\alpha_1,\alpha_2, \alpha_3$ form an inner-complete graph. Then
\begin{eqnarray*}
\alpha_3\alpha_1\alpha_2^k\alpha_3\alpha_1\alpha_2 &=&
\alpha_3\alpha_1\alpha_2^k\alpha_3\lambda\alpha_1 =
\alpha_3\alpha_1\alpha_2^k\mu\alpha_3\alpha_1 \\ &=&
\alpha_3\alpha_1\mu\alpha_2^k\alpha_3\alpha_1 =
\alpha_3\nu\alpha_1\alpha_2^k\alpha_3\alpha_1 =
\alpha_2\alpha_3\alpha_1\alpha_2^k\alpha_3\alpha_1
\end{eqnarray*}
and so $W=W'$.

Let $\alpha_2\alpha_3\alpha_1\alpha_2^k\alpha_3 = \alpha_3 P$.
Then $P \doteq \nu\alpha_1\alpha_2^k\alpha_3$. To remove $\nu$,
$\alpha_1\alpha_2^k\alpha_3$ must start with $\alpha_2$. But the
left canonical form of $\alpha_1\alpha_2^k\alpha_3$ is $(\alpha_1
\alpha_2)(\alpha_2 \alpha_3)\nu^{k-2}$ and so it may start only
with $\alpha_1, \alpha_2,\lambda$. Thus
$\alpha_2\alpha_3\alpha_1\alpha_2^k\alpha_3 \not= \alpha_3
P(\alpha_1\alpha_2\alpha_3)$. Similarly
$\alpha_1\alpha_2^k\alpha_3\alpha_1\alpha_2\not=
Q(\alpha_1\alpha_2\alpha_3)\alpha_1$.

\subsubsection*{\bf II.  Two edges are missing}
In this case, there are two possible graphs as in
Figure~\ref{f:non4} and Figure~\ref{f:non42}. For the graph in
Figure~\ref{f:non4}, let
$W=\alpha_4\alpha_3\alpha_2^{k+1}\alpha_3$ and
$W'=\alpha_1\alpha_2^2\alpha_1^k\alpha_4$.

\putfig{non4}

Add the edges $\lambda,\nu$ so that these edges together with
$\alpha_1,\ldots, \alpha_4$ form an inner-complete graph. Then
\begin{eqnarray*}
\alpha_4\alpha_3\alpha_2^{k+1}\alpha_3 &=&
\lambda\alpha_4\alpha_2^{k+1}\alpha_3 =
\lambda\alpha_2^k\alpha_4\alpha_2\alpha_3 =
\lambda\alpha_2^k\alpha_4\nu\alpha_2 \\ &=&
\lambda\alpha_2^k\alpha_1\alpha_4\alpha_2 =
\lambda\alpha_2^k\alpha_1\alpha_2\alpha_4 =
\lambda\alpha_1\alpha_2\alpha_1^k\alpha_4 =
\alpha_1\alpha_2^2\alpha_1^k\alpha_4
\end{eqnarray*}
and so $W=W'$.

The left canonical form of $\alpha_4\alpha_3\alpha_2^{k+1}$ is
$(\alpha_4\alpha_3)\alpha_2^{k+1}$ and so it can not start with
$\nu$ over the band-generators. Thus
$\alpha_4\alpha_3\alpha_2^{k+1}\not=\alpha_1P(\alpha_1,\ldots,
\alpha_4)$. The right canonical form of
$\alpha_3\alpha_2^{k+1}\alpha_3$ is
$\alpha_3\alpha_2^k(\alpha_2\alpha_3)$ and so it can not finish
with $\alpha_4$. Thus $\alpha_3\alpha_2^{k+1}\alpha_3
\not=Q(\alpha_1,\ldots, \alpha_4)\alpha_4$.

For the graph in Figure~\ref{f:non42}, let $W=\alpha_2 \alpha_3
\alpha_4^k \alpha_2$ and $ W'= \alpha_3 \alpha_4^k \alpha_2
\alpha_3$. Then $W,W'$ satisfy Lemma~\ref{length} by using a
similar argument as for Figure~\ref{f:star1}.

\putfig{non42}

\subsubsection*{\bf III. One edge is missing}
The only possible graph is Figure~\ref{f:non5}. Let
 $W=\alpha_4\alpha_5^k\alpha_2\alpha_3$ and
 $W'=\alpha_1\alpha_4\alpha_5^k\alpha_2$
\putfig{non5}

\begin{eqnarray*} \alpha_4\alpha_5^k\alpha_2\alpha_3 =
\alpha_4\alpha_5^k\lambda\alpha_2 =
\alpha_4\lambda\alpha_5^k\alpha_2 =
\alpha_1\alpha_4\alpha_5^k\alpha_2
\end{eqnarray*}
and so $W=W'$

 The left canonical form of
$\alpha_4\alpha_5^k\alpha_2$ is itself and so it can not start
with $\nu$ over the band-generators. Thus
$\alpha_4\alpha_5^k\alpha_2\not=\alpha_1P(\alpha_1,\ldots,
\alpha_5)$. Let $\alpha_5^k\alpha_2\alpha_3\doteq Q\alpha_2$ then
$Q\doteq\alpha_5^k\nu$ and $\alpha_5,\nu$ commute, so $\nu$ can
not be removed. Thus $\alpha_5^k\alpha_2\alpha_3
\not=Q(\alpha_1,\ldots, \alpha_5)\alpha_2$.

\subsubsection*{\bf IV. No edges are missing}

If we add any edge to the inner-complete graph, multiple edges
are created. Thus only the inner-complete graph is eligible.

\end{proof}

\begin{thm}
Among all linearly spanned graphs without multiple edges, only the
Artin graphs and the inner-complete graphs have the embedding
property.
\end{thm}

\begin{proof}
We have already discussed about linearly spanned graphs with 3
vertices. Let $\Gamma$ be a linearly spanned graph with more than
3 vertices that has the embedding property. First we choose a
connect full subgraph $\Gamma'$ with 4 vertices from $\Gamma$ so
that there is a separating circle satisfying the hypothesis of
Theorem~\ref{circle}. Choose 4 vertices that form a connected
subtree in a spanning tree of $\Gamma$. Take the full subgraph
with these 4 vertices. If there is no other vertices in faces of
this full subgraph, then this is a desired full subgraph. If there
is other vertices in the faces of this full subgraph and none of
them is adjacent to the chosen 4 vertices, then there is an
edge-path starting at a vertex $v$ on a face and ending at one of
the 4 vertices such that the edge-path intersect the full subgraph
of the 4 vertices and so $v$ is contained in a pseudo face. Since
$\Gamma$ is linearly spanned, this can not happen. Thus at least
one of vertices on faces, say $w$, is adjacent to one of the 4
vertices. Then we have the less number of unwanted vertices
contained in faces of a new connected full subgraph that is
obtained by replacing one of the 4 vertices by $w$. By repeating
this process, we eventually obtain a connect full subgraph
$\Gamma'$ with 4 vertices $v_1,v_2,v_3,v_4$ that can be separated
by a circle from other vertices of $\Gamma$. Theorem~\ref{circle}
and Theorem~\ref{4vertex} say that $\Gamma'$ must be either the
Artin graph or the inner-complete graph with 4 vertices. Choose
a vertex $v_5$ in $\Gamma-\Gamma'$ that is adjacent to one of
$v_1,v_2,v_3,v_4$. If $\Gamma'$ is an inner-complete graph, then
each full subgraph with the 4 vertices consisted of $v_5$ and any
three vertices from $v_1,v_2,v_3,v_4$ must be an inner-complete
graph by Theorem~\ref{circle} and Theorem~\ref{4vertex}. Thus the
full subgraph with the 5 vertices $v_1,v_2,v_3,v_4,v_5$ is an
inner-complete graph. By repeating this process, $\Gamma$ itself
eventually becomes an inner-complete graph.

We now suppose that $\Gamma'$ is the Artin graph with 4 vertices
$v_1,v_2,v_3,v_4$ where $v_1,v_4$ are of valency 1 and the other
vertices are of valency 2.  Choose a vertex $v_5$ in $\Gamma-\Gamma'$ that
is adjacent to any one of $v_1,v_2,v_3,v_4$. By Theorem~\ref{circle}
and Theorem~\ref{4vertex}, $v_5$ can be adjacent to either one or both
of vertices $v_1,v_4$. If $v_5$ is adjacent to both vertices, we stop
the process. If $v_5$ is adjacent to one vertex, say $v_4$,
then $v_1,\ldots,v_5$ form an Artin graph and we repeat this process.
Eventually we see that $\Gamma$ must contain one of the following
two types of graphs as a full subgraph separated by a circle
as in Theorem~\ref{circle} unless $\Gamma$ itself is an Artin graph.
The proof will be completed when we show
that both types of graphs do not have the embedding property.

\subsection*{\boldmath $m$-Gon} The graph of this type is depicted in
Figure~\ref{f:kcycle}. Let $W=\alpha_1 \alpha_2 \cdots \alpha_i^k
\cdots \alpha_{m-1}$ and $W'= \alpha_2 \cdots \alpha_i^k \cdots
\alpha_m$. Then one can show that $W,W'$ satisfy the hypothesis of
Lemma~\ref{length} by a similar but longer argument as for
Figure~\ref{f:cycle}. Thus $\Gamma$ does not have the embedding
property

 \putfig{kcycle}

\subsection*{\boldmath Pseudo $m$-gon} The graph of this type
is depicted in Figure~\ref{f:nonpl2}. Let $$W=\alpha_2 \alpha_3 \cdots
\alpha_m \alpha_1 \alpha_2^k \alpha_3 \cdots \alpha_{m-1} \alpha_m
\alpha_1 \alpha_2 \cdots \alpha_{m-2}$$ and $$W'= \alpha_3
\alpha_4 \cdots \alpha_m \alpha_1 \alpha_2^k \alpha_3 \cdots
\alpha_{m-1} \alpha_m \alpha_1 \alpha_2 \cdots \alpha_{m-1}$$

\putfig{nonpl2}

Add supplementary edges $\lambda, \mu, \nu. \alpha'_{m-1},
\tau_{m-1}, \tau_m$ as in Figure~\ref{f:nonpl2} to form a subgraph
of the inner-complete graph with $m+1$ vertices. Then
\begin{eqnarray*} \lefteqn{\alpha_2 \alpha_3
\cdots \alpha_m \alpha_1 \alpha_2^k \alpha_3 \cdots \alpha_{m-1}
\alpha_m \alpha_1 \alpha_2 \cdots \alpha_{m-2}}
\\ &=& \alpha_3 \alpha_4 \cdots \alpha_m \lambda \alpha_1
\alpha_2^k \alpha_3 \cdots \alpha_{m-1} \alpha_m \alpha_1 \alpha_2
\cdots \alpha_{m-2}\\ &=& \alpha_3 \alpha_4 \cdots \alpha_m
\alpha_1 \mu \alpha_2^k \alpha_3 \cdots \alpha_{m-1} \alpha_m
\alpha_1 \alpha_2 \cdots \alpha_{m-2} \\ &=& \alpha_3 \alpha_4
\cdots \alpha_m \alpha_1 \alpha_2^k \alpha_3 \cdots \alpha_{m-1}
\mu \alpha_m \alpha_1 \alpha_2 \cdots \alpha_{m-2}\\ &=& \alpha_3
\alpha_4 \cdots \alpha_m \alpha_1 \alpha_2^k \alpha_3 \cdots
\alpha_{m-1} \alpha_m \nu \alpha_1 \alpha_2 \cdots \alpha_{m-2}\\
&=& \alpha_3 \alpha_4 \cdots \alpha_m \alpha_1 \alpha_2^k \alpha_3
\cdots \alpha_{m-1} \alpha_m \alpha_1 \alpha_2 \cdots \alpha_{m-1}
\end{eqnarray*}
and so $W=W'$.

We now show that $$\alpha_2 \alpha_3 \cdots \alpha_m \alpha_1
\alpha_2^k \alpha_3 \cdots \alpha_{m-1} \alpha_m \alpha_1 \alpha_2
\cdots \alpha_{m-3}\not=\alpha_3P(\alpha_1,\ldots,\alpha_m).$$
Suppose that
$$\alpha_2 \alpha_3 \cdots \alpha_m \alpha_1 \alpha_2^k
\alpha_3 \cdots \alpha_{m-1} \alpha_m \alpha_1 \alpha_2 \cdots
\alpha_{m-3}=\alpha_3P.$$ Then  \begin{eqnarray*} P &\doteq&
\alpha_4 \cdots \alpha_m\lambda \alpha_1\alpha_2^k \alpha_3 \cdots
\alpha_{m-1} \alpha_m\alpha_1\alpha_2\cdots\alpha_{m-3}\\ &\doteq
& U_m\tau_m
\end{eqnarray*}
where $U_m\equiv\alpha_4 \cdots \alpha_m\alpha_1\alpha_2^k
\alpha_3 \cdots \alpha_{m-1}
\alpha_m\alpha_1\alpha_2\cdots\alpha_{m-3}$. Due to the triangular
relation $\alpha_{m-1}\tau_m=\alpha_{m-2}\alpha_{m-1}$, the
positive word $U_m$ must end with $\alpha_{m-1}$ in order for
$\tau_m$ to disappear. We will show that $U_m$ can not end with
$\alpha_{m-1}$ by an induction on $m\ge 3$. Clearly
$U_3=\alpha_1\alpha_2^k \alpha_3$ can not end with $\alpha_2$ in
the band-generator presentation.

\begin{eqnarray*}
U_m &\doteq& \alpha_4 \cdots
\alpha_{m-1}\alpha_m\alpha_1\alpha_2^k \alpha_3 \cdots
\alpha_{m-1} \alpha_m\alpha_1\alpha_2\cdots\alpha_{m-3}\\
 &\doteq&
\alpha_4\cdots\alpha_{m-2}\alpha'_{m-1}\alpha_{m-1}\alpha_1
\alpha_2^k \alpha_3 \cdots \alpha_{m-1}
\alpha_m\alpha_1\alpha_2\cdots\alpha_{m-3}\\ &\doteq&
\alpha_4\cdots\alpha_{m-2}\alpha'_{m-1}\alpha_1 \alpha_{m-1}
\alpha_2^k \alpha_3 \cdots \alpha_{m-1}
\alpha_m\alpha_1\alpha_2\cdots\alpha_{m-3}\\ &\doteq&
\alpha_4\cdots\alpha_{m-2}\alpha'_{m-1}\alpha_1 \alpha_2^k\alpha_3
\cdots \alpha_{m-1}\alpha_{m-2}\alpha_{m-1}
\alpha_m\alpha_1\alpha_2\cdots\alpha_{m-3}\\ &\doteq&
\alpha_4\cdots\alpha_{m-2}\alpha'_{m-1}\alpha_1 \alpha_2^k\alpha_3
\cdots \alpha_{m-2}\alpha_{m-1}\alpha_{m-2}
\alpha_m\alpha_1\alpha_2\cdots\alpha_{m-3}\\ &\doteq&
\alpha_4\cdots\alpha_{m-2}\alpha'_{m-1}\alpha_1
\alpha_2^k\alpha_3\cdots \alpha_{m-2}
\alpha_{m-1}\alpha_m\alpha_1\alpha_2\cdots\alpha_{m-4}\alpha_{m-2}\alpha_{m-3}\\&\doteq&
\alpha_4\cdots\alpha_{m-2}\alpha'_{m-1}\alpha_1
\alpha_2^k\alpha_3\cdots \alpha_{m-2}
\alpha'_{m-1}\alpha_{m-1}\alpha_1\alpha_2\cdots\alpha_{m-4}\alpha_{m-2}\alpha_{m-3}\\&\doteq&
\alpha_4\cdots\alpha_{m-2}\alpha'_{m-1}\alpha_1
\alpha_2^k\alpha_3\cdots \alpha_{m-2}
\alpha'_{m-1}\alpha_1\alpha_2\cdots\alpha_{m-4}\alpha_{m-1}\alpha_{m-2}\alpha_{m-3}\\&\doteq&
U_{m-1}\alpha_{m-1}\alpha_{m-2}\alpha_{m-3}
\end{eqnarray*}
It is easy to check that
$U_{m-1}\alpha_{m-1}\alpha_{m-2}\alpha_{m-3}$ ends with
$\alpha_{m-1}$ if and only if $U_{m-1}\alpha_{m-1}\alpha_{m-2}$
ends with $\alpha_{m-1}$ if and only if $U_{m-1}$ ends with
$\alpha_{m-2}$ in the band-generator presentation. By the
induction hypothesis, $U_{m-1}$ can not end with $\alpha_{m-2}$.
Similarly we have that $$ \alpha_3 \cdots \alpha_m \alpha_1
\alpha_2^k \alpha_3 \cdots \alpha_{m-1} \alpha_m \alpha_1 \alpha_2
\cdots \alpha_{m-2}\not=Q(\alpha_1,\ldots,\alpha_m)\alpha_{m-1}.$$
Consequently $W,W$ satisfy the hypothesis of Lemma~\ref{length}.
\end{proof}

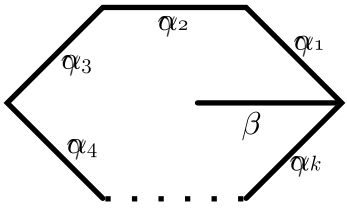
\end{document}